%Paper: alg-geom/9406001
%From: Yukari-i <yukari-i@tansei.cc.u-tokyo.ac.jp>
%Date: Sat, 11 Jun 94 17:06:26 JST
%Date (revised): Fri, 22 Jul 94 17:21:48 JST

%
\input amstex
\mag=\magstep1
\documentstyle{amsppt}
\NoRunningHeads
\vsize=23.5true cm
\hsize=15.5true cm
\topmatter
\title Gorenstein Quotient Singularities of Monomial Type in Dimension Three
\endtitle
\author Yukari ITO \endauthor
\affil Department of Mathematical Sciences \\
 University of Tokyo \\
 7-3-1 Hongo, Bunkyo, Tokyo, 113, Japan \\
 yukari-i{\@}tansei.cc.u-tokyo.ac.jp \endaffil
\endtopmatter
\document
\subhead  \nofrills{\bf
\S 1. Introduction}
\endsubhead \par
The purpose of this paper is to construct a crepant resolution of quotient
singularities by finite subgroups of $SL(3,\Bbb C)$ of monomial type
(type (B),(C) and (D) in [5]), and prove that the Euler number of the
resolution is equal to the number of conjugacy classes.
 This preprint contains the results which I talked at Research Institute for
Mathematical Sciences of Kyoto University on 13th May, 1994.
\par
 The problem of finding a nice resolution of quotient singularities by finite
subgroups of $SL(3,\Bbb C)$ arose from mathematical physics.
\definition{Definition}(Orbifold Eular characteristic) In the physics of
superstring theory, one considers the string propagation on a
manifold $M$ which is a quotient by a finite subgroup of symmetries $G$. By a
physical argument of string vacua of $M/G$, one concludes that the correct
Euler number for the theory should be
the ``orbifold Euler characteristic"[2,3], defined by \par
$$ \chi(M,G) = \frac 1{\vert G \vert} \sum_{gh=hg} \chi(M^{<g,h>}),$$
 where the summation runs over all pairs of commuting elements of $G$,
and $M^{<g,h>}$ denotes the common fixed set of $g$ and $h$. For the physcist's
interest, we only consider $M$ whose quotient space $M/G$ has trivial canonical
\enddefinition
\proclaim{ Conjecture I} {\rm([2,3])}\par
\it{There exists a resolution of singularities $\widetilde{M/G}$ s.t.
$\omega_{\widetilde{M/G}} \simeq \Cal O_{\widetilde{M/G}}$, and
$$ \chi(\widetilde{M/G}) = \chi(M,G). $$} \endproclaim
This conjecture follows from its local form [4]:

\proclaim{ Conjecture II} {\rm(local form)}\par
\it{Let $G\subset SL(3,\Bbb C)$ be a finite group. Then there exists a
resolution
 of singularities\par $\sigma$ : $\widetilde X \longrightarrow {\Bbb C}^3/G$
with
$\omega_{\widetilde X} \simeq \Cal O_{\widetilde X} $ and
$$ \chi(\widetilde X) = \sharp \{ \text{conjugacy class of } G \}. $$
} \endproclaim
In algebraic geometry, the conjecture says that a minimal model of the quotient
space by a finite subgroup of $SL(3,\Bbb C)$ is non-singular.\par

Conjecture II  was proved for abelian groups  by Roan ([10]), and
independently by Markushevich, Olshanetsky and Perelomov ([8]) by using toric
method. It was also proved for 5 other groups, for which $X$ are hypersurfaces:
(i) $W{A_3}^+, W{B_3}^+, W{C_3}^+,$ where $WX^+$ denotes the positive
determinant part of the Weyl group $WX$ of a root system $X$ by Bertin and
Markushevich
 ([1]), (ii) $H_{168}$  by Markushevich ([7]), and
(iii) $I_{60}$  by Roan ([11]).
Recently I proved Conjecture II for trihedral groups [5]:
\definition{ Definition} {\it A trihedral group} is a finite group
 $G=<H,T>\ \subset SL(3,\Bbb C)$,
where $H \subset SL(3,\Bbb C)$ is a finite group generated by diagonal matrices
and
$$ T=\pmatrix
     0 & 1 & 0 \\
     0 & 0 & 1 \\
     1 & 0 & 0 \\
    \endpmatrix. $$
\enddefinition
\definition{Definition}
{\it Trihedral singularities} are quotient singularities by trihedral groups.
\enddefinition
\definition{Definition}
A resolution of singularities $f:Y \longrightarrow X$ of a normal variety $X$
with $K_X$ being $\Bbb Q$-Cartier is {\it crepant} if $K_Y=f^*K_X$.
\enddefinition
\proclaim{Theorem 1.1[5] }\par
\it{Let $X=\Bbb C^3/G$ be a quotient space by a trihedral group $G$. Then there
exists a crepant resolution of singularities
$$f:\widetilde X \longrightarrow X,$$
 and
$$ \chi(\widetilde X) = \sharp \{ \text{conjugacy class of }G \}. $$}
\endproclaim
\proclaim{Theorem 1.2 (Main Theorem)}\par
\it{The conjecture II holds for the following groups: \par
I. $G_1 = < H, S>$ \par
II. $G_2 = < H, H', S> $ \par
III. $G_3 = < H, S, T >$  $(r\not\equiv0$(mod 3))\par
IV. $G_4= <G_3, C>$ \par
V. $G_5=<C,S>$ \par
where $ H = \frac 1{r}( 0,1,-1), H' = \frac 1{r}( 1,-1,0), C=\frac1{3}(1,1,1)$
and
$$ S =\pmatrix
       -1 & 0 & 0 \\
        0 & 0 & -1 \\
        0 & -1 & 0 \\
       \endpmatrix .$$}
\endproclaim
\remark{Remark 1.3}
 These singularities are different from trihedral singularities, but main idea
for these proofs are based on the method of trihedral case.
If we take any abelian group for $H$ in Theorem 1.2, we can use the same method
 and the conjecture II is true.
\endremark

In this paper, we prove Theorem 1.2.I  in section 2,  II in section 3, III in
section 4 and IV and V in section 5.

\subhead \nofrills{\bf
\S 2. Proof of I
}\endsubhead
\par

 \proclaim {Proposition 2.1} \par
\it{Let $X = {\Bbb C} ^3/ G$, and $Y = {\Bbb C}^3 /G'$. Then there exists the
following diagram: }
$$ \CD
   @. @. {\widetilde X} \\
   @. @. @VV{\tau}V     \\
   @. {\widetilde Y} @>{\widetilde\mu}>> {\widetilde Y}/{\Bbb Z_2} \\
   @. @VV{\pi}V @VV{\widetilde\pi}V \\
   {\Bbb C^3} @>>> {\Bbb C^3}/{G'}=Y @>{\mu}>> {\Bbb C^3}/G=X
 \endCD $$
\it{where $\pi$ is a resolution of the singularity of $Y$, and $\widetilde\pi$
is the induced morphism, $\tau$ is a resolution of  the singularity by $\Bbb
Z_2$, and $\tau \circ \widetilde \pi$ is a crepant resolution of the
singularity of $X$.}\endproclaim

\demo{\bf{Sketch of the proof}} As a resolution $\pi$ of $Y$, we take a toric
resolution, which is also crepant. Then we lift the $\Bbb Z_2$-action on $Y$ to
its crepant resolution $\widetilde Y$ and form the quotient $\widetilde Y/\Bbb
Z_2$. This quotient gives in a natural way a partial resolution of the
singularities of $X$. The crepant resolution $\widetilde X \longrightarrow
\widetilde Y/\Bbb Z_2$ of the singularities
of $\widetilde Y/\Bbb Z_2$ induces a complete resolution of $X$.\par
Under the action of  $\Bbb Z_2$, the singularities of $\widetilde Y/\Bbb Z_2$
lie in the union of the image of the exceptional divisor of $\widetilde Y$
under $\widetilde Y  \longrightarrow \widetilde Y/\Bbb Z_2$ and the locus
$C:(x=0,y=-z)$ in $\Bbb C^3$. \par
In the resolution $\widetilde Y$ of $Y$, the group $\Bbb Z_2$ permutes
exceptional divisors. So the fixed points on the exceptional divisors consist
of one point or three points. \enddemo

Now, we see the proof more precisely.
\par
At first, let $G'$ be an abelian group which generated by diagonal matrices in
$G$. It is a normal subgroup of $G$ in any case of Theorem 1.2.

Then we recall a toric resolution when $G'$ is a finite abelian group in
$SL(3,\Bbb C)$.\par
Let $\Bbb R^3$ be a vector space, and $\{ e^i | i=1,2,3\}$ the standard base.
For all $v=\frac 1{r}(a,b,c)\in G'$, $L= $the lattice generated by $e^1,e^2$
and $e^3$,  $N:=L+\sum \Bbb Zv$,
$$ \sigma = \left\{ \sum_{i=1}^3 x_ie^i \in \Bbb R^3, \quad x_i \geq 0, \forall
i \right\}. $$
We regard $\sigma$ as a rational convex polyhedral cones in $N_{\Bbb R}$. The
corresponding affine torus embedding $X_{\sigma}$ is defined as Spec$\Bbb
C[\check{\sigma} \cap M]$, where $M$ is the dual lattice of $N$ and
$\check{\sigma}$ is a dual cone of $\sigma$ in $M_{\Bbb R}$ defined by
$\check{\sigma} =\{ \xi\in M_f{\Bbb R} | \xi(x) \geq 0 , \forall x \in \sigma\}
$.  \par
$\qquad \qquad \varDelta :=$ the simplex in $N_{\Bbb R}$
$$ \left\{\sum_{i=1}^3 x_ie^i  \quad ; x_i \geq 0,\quad \sum_{i=1}^3 x_i=1
\right\} \qquad \qquad$$ \par
$$  t:N_{\Bbb R} \longrightarrow \Bbb R \qquad \sum_{i=1}^3x_ie^i \longmapsto
\sum_{i=1}^3 x_i \ \qquad \qquad \qquad  $$\par
$$ \Phi := \left\{ \frac 1{r}(a,b,c)\in G'\  | \ a+b+c=r \right\}
\qquad \qquad \qquad $$
\proclaim{Lemma 2.2} \par
\it{$Y=\Bbb C^3/G'$ corresponds to $\sigma$ in $N=L+ \sum_{v \in \Phi}\Bbb
Zv$.}
\endproclaim
\demo{\bf{Proof}}
Since $Y=\text{Spec}(\Bbb C[x,y,z]^{G'})$, $x^iy^jz^k$ is $G'$-invariant if and
only if $ \alpha i+\beta j + \gamma k \in \Bbb Z$ for all
$(\alpha,\beta,\gamma) \in G'$.  \qed \enddemo
\example{Remark 2.3}
Let $\frac 1{r}(a,b,c) \not= (0,0,0)$ be an element of $G'$. There are two
types.\roster
\item \quad $abc \not= 0$
\item \quad $abc = 0$
\endroster
We denote by $G_1$ (resp. $G_2$) the set of the elements of type (1) (resp.
(2)). So $G'\backslash \{e\} = G_1 \amalg G_2$. \par
There are also two types in $\Phi$ as above. So we denote by $\Phi_1$ (resp.
$\Phi_2)$ the set of lattice points of type (1) (resp. (2)). Then $\Phi=\Phi_1
\amalg \Phi_2$.\par
Let $\lambda_i$ be maps from $G_i$ to $\Phi_i$ (i=1,2).
$$ \lambda_i\  :\ G_i \ \longrightarrow \ \Phi_i $$
where $\lambda_1$ maps $g=\frac 1{r}(a,b,c)\ (a+b+c=r)$ and $g^{-1}= \frac
1{r}(r-a,r-b,r-c)$ to a lattice point $\frac 1{r}(a,b,c)$, and $\lambda_2$ maps
$g=\frac 1{r}(a,b,c)$ to a lattice point $\frac 1{r}(a,b,c)$.
 $$ \left\{G_1\right\} \overset {2:1} \to \longrightarrow
\left\{\Phi_1\right\}\ , \qquad \left\{G_2\right\} \overset {1:1} \to
\longrightarrow \left\{\Phi_2\right\}. $$
 \par
Therefore there exist a correspondence between the sets of elements of
$G'-\{e\}$ and $\Phi$, which is 2:1 on $G_1$ and 1:1 on $G_2$. $\Phi$
corresponds to the  exceptional divisors of a toric resolution given below.
\endexample
\proclaim{Claim I} \it{There exists a toric resolution of $Y$ where $\Bbb Z_2$
acts symmetrically on the exceptional divisors.}
\endproclaim
\demo{\bf{Proof}}  We can construct a unique simplicial decomposition.\enddemo
\proclaim{Claim II}
\it{Let $X_S$ be the corresponding torus embedding, then $X_S$ is non-singular.
}\endproclaim
\demo{\bf{Proof}}
It is sufficient to show that the $\sigma(s)$ are basic.
Let $ w^1,w^2,w^3 \in \Phi \cup \{e^i\}_{i=1}^3$ which are linearly independent
over $\Bbb R$. Assume that the simplex
$$ \left \{ \sum_{i=1}^3 \alpha_iw^i \ | \ \alpha_i \geq 0 ,\ \sum_{i=1}^3
\alpha_i=1 \right \} $$
intersects $\Phi \cup \{e^i\}_{i=1}^3$ only at $\{w^i\}_{i=1}^3.$
\par
Lattice $N_0$ generated by $\{w^i\}_{i=1}^3$ is sublattice of $N$.
If we assume $N\not= N_0$, then there exists
$\beta=\beta_1w^1+\beta_2w^2+\beta_3w^3 \ \in\  N\backslash N_0$\  $( 0 \leq
\beta_i  < 1,\  \beta_i \in \Bbb R$ and the strict inequality holds at least
for one $i$. ) \par
$t(\beta)=\sum\beta_i t(w^i) = \sum \beta_i$, $0 < \sum \beta_i < 3$ and $t(N)
\in \Bbb Z$, then $t(\beta)=1$ or 2.
If $t(\beta)=2$, then we can replace it by $\beta'=
\sum_{i=1}^3(1-\beta_i)w^i$,  so we can assume that $t(\beta)=1$. \par
Now, there exists an element  $\beta$ in $
\{ \sum \alpha_iw^i \ | \ \alpha_i \geq 0, \ \sum \alpha_i =1\} \cap (N -
N_0)$,
 which is contained in $\varDelta \cap N$.
Since $$N = \left\{\bigcup \Sb v \in \Phi \endSb ( v\oplus L)\right\}\bigcup L,
$$
$\varDelta \cap N = \Phi \cup \{e^i\}_{i=1}^3.$ From our assumption, we
conclude that $\beta=w^i$ for some $i$, which contradicts $\beta\not \in N_0$.
Therefore $N=N_0$.
\qed \enddemo
We obtain a crepant resolution $\pi_S:X_S\longrightarrow \Bbb C^3/G'$, because
$X_S$ is non-singular.

\proclaim{Claim III}
\it{Let $F$ be the fixed locus on $\widetilde Y$ under the action of $\Bbb
Z_2$, and $E$ be exceptional divisors of $\widetilde Y \longrightarrow Y$. Then
$$ F_0 := F \cap E = \cases \text{1 point \qquad \ : $G'$ is type (I)} \\
                                 \text{2 points \qquad : $G'$ is type (II)}.
\endcases $$}
\endproclaim
\demo {\bf{Proof}}
Considering the dual graph of exceptional divisors by the toric resolution and
from Remark 2.3, we can identify the two exceptional divisors by the action of
$\Bbb Z_2$ except the central component which is a component in the center
of the exceptional locus. Then there are 2 possibilities of the central
component; \par
 Type (I): one point.\par
 Type (II): a divisor which is isomorphic to $\Bbb P^1$.
Let $(y:z)$ be a coordinate of $\Bbb P^1$, then the action of $\Bbb Z_2$ is
$$ (y:z) \longmapsto (-z:-y). $$
Then the number of the fixed points are two, whose coordinates are
$$(1:1),(1:-1)\qed$$
\enddemo
Furthermore the $\Bbb Z_2$-action in the neighbourhood of a fixed point is
analytically isomorphic to some linear action. \par
Now, we consider the resolution of the singularity of $\Bbb C^3$ by the group
$\Bbb Z_2$.\par
$F_0$ consists of 1 or 2 points, and $F$ = $F_0 \cup C'$ where $C'$ is a strict
transform of the fixed locus $C$ in $Y$ under the action of $\Bbb Z_2$. \par

\proclaim {Claim IV} \it{Let $Z=\Bbb C^3/\Bbb Z_2$, then $\chi(\widetilde Z) =
\chi( \Bbb C^3, \Bbb Z_2) = 2.$}
\endproclaim
\demo{\bf{Proof}}
There is a representation of $\Bbb Z_2$ in $SL(3,\Bbb C)$:
$$S' = \pmatrix
         1 & 0 & 0 \\
         0 & -1 & 0 \\
         0 & 0 & -1 \\
        \endpmatrix $$
The quotient singularities by $S'$ are $A_1 \times \{x-$axis\} which are not
isolated and the exceptional divisor
of the resolution is $\Bbb P^1$-bundle.
\qed \enddemo

\proclaim {Claim V} \it{
The resolution  $\tau \circ \widetilde \pi $ is a crepant resolution.}
\endproclaim

\proclaim {Lemma 2.4} \par
\it{Let $X:= \Bbb C^3/<G',T>$, and $f:\widetilde X \longrightarrow X$ the
crepant resolution as above. Then the Euler number of $\widetilde X$ is given
by
$$ \chi(\widetilde X) = \frac 1{2} (\vert G' \vert -k)+2k $$
where $$ k  =
       \cases
           1 \qquad \text{if\ \   $|G'| \equiv 1$ (mod 2) \  (type (I))}  \\
           2 \qquad \text{if\ \   $|G'| \equiv 0$ (mod 2) \  (type (II))}.
       \endcases $$}
\endproclaim

\demo{\bf{Proof}}
For an abelian group $G'$, we have a toric resolution
$$ \pi : \widetilde Y \longrightarrow Y = \Bbb C^3/G', $$
and $\chi(\widetilde Y) = \vert G'\vert$ ([8],[11]). \par
By the action of $\Bbb Z_2$, the number of fixed points in the exceptional
divisor by $\sigma$ is equal to $k$, hence
$$ \chi(\widetilde Y/\frak A_2)=\frac 1{2} (\vert G' \vert -k)+k. $$
By the resolution of the fixed loci, Eular characteristic of the each
exceptional locus is 2. (Claim IV) \par
Therefore,
$$ \chi(\widetilde X) = \frac 1{2}(\vert G'\vert -k) + 2k.  \qed$$
\enddemo

\proclaim{Theorem 2.5} \par
\it{$$\chi (\widetilde X) = \sharp \{ \text{conjugacy class of $G$} \}. $$
}\endproclaim
\demo{\bf{Proof}}
\roster
 \item  Case (I) : $\vert G'\vert = 2m+1 , (0 <m\in \Bbb Z)$
\par For a nontrivial element $g \in G'$,
there are two conjugate elements
$g$ and   $SgS$.
There are $m$ couples of this type.
There are 2 other  conjugacy classes $[e]$ and  $ [S]$.
 Therefore, there are $m+2$ conjugacy classes in $G$. \par
Then
$$ \allowdisplaybreaks \align
   \chi(\widetilde X) &= \frac 1{2}(|G'|-1)+2 \\
                      &= m+2 \\
                      &= \sharp \{ \text{ conjugacy class of } G \}.
\endalign $$
 \item  Case (II) : $\vert G'\vert = 2m , (0 <m \in \Bbb Z)$
\par There are 2 elements in the center of $G'$: $e, a=\frac 1{2}(0,1,1)$.
 The other  $2m-2$ elements in $G'$ are divived into $m-1$ conjugacy classes as
in (1). There are 4 other conjugacy classes
$e,a,[S],[aS]$.
Therefore, there are $m+3$ conjugacy classes in $G$. \par
Then
$$ \allowdisplaybreaks \align
     \chi(\widetilde X) & =\frac 1{2} (\vert G'\vert -2) +4 \\
                        & = m+3 \\
                        & = \sharp \{ \text{ conjugacy class of } G \}. \qed
\endalign     $$
\endroster
\enddemo

\subhead \nofrills {\bf \S 3.Proof of II } \endsubhead \par
In this case, we do similarly as I. So we check the toric resolution:
\proclaim{Proposition 3.1} \it{There exists a toric resolution of $Y$ where
$\Bbb Z_2$ acts symmetrically on the exceptional divisors.}
\endproclaim
\demo{\bf{Proof}}  We show that we can construct a simplicial decomposition
$\{ \sigma(s)\}_{s \in S}$ of the simplex $\varDelta$ which is $\Bbb
Z_2$-invariant and the set of  vertices is exactly $\Phi \cup \{e^i\}_{i=1}^3$.
\par
 Let us consider the distance $d$ between $\frac 1{r}(a,b,c)$ and $\frac
1{2}(a,h,h)$ $(h=(b+c)/2)$ given by
$$ d\left(\frac 1{r}(a,b,c)\right)=\left\vert\frac b{r}-\frac 1{h} \right\vert
                        + \left\vert \frac c{r}- \frac 1{h}\right\vert.$$
\roster
\item Find lattice points $P_i= \frac 1{r}(a,b,c)$ whose distance $d$ is the
minimum among the points for each $a$ in $\Phi$ in the domain $D=\{1/2\geq
y$\}.
\item Make a triangle whose vertices are $P_i$ and  ${P_i}'= \frac 1{r}(a,c,b)$
symmetrically.
\item Decompose $D$ whose vertices are $P_i$,\ ${P_i}'$,\ $(1,0,0)$  and
$(0,0,1)$, into  simplexes using vertices in $\Phi$.
 We call this decomposition $S_1$.
\item By the action of $\Bbb Z_2$, we obtain $S_2$ on the other
triangle. Therefore we obtain a ``symmetric"
resolution.\qed
\endroster
 \enddemo
Next, we see the singularities in $\widetilde Y/\Bbb Z_2$:

\proclaim{Proposition 3.2}
\it{Let $F$ be the fixed locus on $\widetilde Y$ under the action of $\Bbb
Z_2$, and $E$ be exceptional divisors of $\widetilde Y \longrightarrow Y$. Then
$$ F_0 := F \cap E = \text{\{r points \}} $$}
\endproclaim
\demo {\bf{Proof}}
Considering the dual graph of exceptional divisors by the toric resolution and
from Remark 2.3, we can identify the two exceptional divisors by the action of
$\Bbb Z_2$ except the central components which are near central part ( from
(1,0,0)
 to $\frac 1{2}(0,1,1)$). Then there are 2 possibilities of the central
component locally; \par
 Type (I): one point.\par
 Type (II): a divisor which is isomorphic to $\Bbb P^1$.
{}From analogue of Claim III in section 2, there are $n$ points in $F_0.$
\enddemo
Furthermore the $\Bbb Z_2$-action in the neighbourhood of a fixed point is
analytically isomorphic to some linear action. \par
\proclaim {Lemma 3.3} \par
\it{Let $X:= \Bbb C^3/<G',T>$, and $f:\widetilde X \longrightarrow X$ the
crepant resolution as above. Then the Euler number of $\widetilde X$ is given
by
$$ \chi(\widetilde X) = \frac 1{2} (\vert G' \vert -r)+2r $$
\endproclaim

\proclaim{Theorem 3.4} \par
\it{$$\chi (\widetilde X) = \sharp \{ \text{conjugacy class of $G$} \}. $$
}\endproclaim
\demo{\bf{Proof}}

$\vert G'\vert = r^2 , (0 <m\in \Bbb Z)$
\par
There are $r$ elements of type $\frac 1{r}(a_i,h_i,h_i)=:a_i$. So
for other nontrivial element $g \in G'$,
there are two conjugate elements
$g$ and $SgS$.
There are $\frac {r^2-r}{2}$ couples of this type.
There are $r$ other  conjugacy classes $ [S]$ and  $[a_iS]$.
 Therefore, there are $\frac{r^2-r}{2}+2r$ conjugacy classes in $G$. \par
Then
$$   \chi(\widetilde X) = \sharp \{ \text{ conjugacy class of } G \}.$$
\enddemo

\subhead \nofrills {\bf \S 4. Proof of III  } \endsubhead \par
In this section, we assume that $r \equiv 1$ or $2$ (mod 3).
\proclaim{Proposition 4.1} \par
\it{Let $X = {\Bbb C} ^3/ G$, and $Y = {\Bbb C}^3 /G'$. Then there exists the
following diagram: }
$$ \CD
   @. @. {\widetilde X} \\
   @. @. @VV{\tau}V     \\
   @. {\widetilde Y} @>{\widetilde\mu}>> {\widetilde Y}/{\frak S_3} \\
   @. @VV{\pi}V @VV{\widetilde\pi}V \\
   {\Bbb C^3} @>>> {\Bbb C^3}/{G'}=Y @>{\mu}>> {\Bbb C^3}/G=X
 \endCD $$
\it{where $\pi$ is a resolution of the singularity of $Y$, and $\widetilde\pi$
is the induced morphism, $\tau$ is a resolution of  the singularity by $\frak
S_3$, and $\tau \circ \widetilde \pi$ is a crepant resolution of the
singularity of $X$.}\endproclaim
We do similarly as in \S 2,3.
\proclaim{Claim I} \it{There exists a toric resolution of $Y$ where $\frak S_3$
acts symmetrically on the exceptional divisors.}
\endproclaim
\demo{\bf{Proof}}  We construct a simplicial decomposition
 with conditions for trihedral case [5] and Proposition 3.1 in \S 3.\enddemo
\proclaim{Claim II}
\it{Let $F$ be the fixed locus on $\widetilde Y$ under the action of $\frak
S_3$, and $E$ be exceptional divisors of $\widetilde Y \longrightarrow Y$. Then
$$ F_0 := F \cap E = \text{\{ 3r-2 points \}},
 $$
where one of them are by $\frak S_3$, the others are by $\Bbb Z_2$.}
\endproclaim
\demo {\bf{Proof}}
Considering the dual graph of exceptional divisors by the toric resolution and
from Remark 2.3, we can identify the two exceptional divisors by the action of
$\Bbb Z_2$ except the central component which is a component in the center
of the exceptional locus. And the central component is a point, which fixed by
the action of $\frak S_3.$
\enddemo
\proclaim {Claim III} \it
{Let $Z=\Bbb C^3/\frak S_3$, then $\chi(\widetilde Z) = \chi( \Bbb C^3, \frak
S_3) = 3.$}
\endproclaim
\demo{\bf{Proof}}
There is a equivalent representation of $T$ in $SL(3,\Bbb C)$:
$$T' = \pmatrix
         1 & 0 & 0 \\
         0 & \omega & 0 \\
         0 & 0 & \omega^2 \\
        \endpmatrix $$
The quotient singularities by $<S,T'>$ are same as the case of Theorem 1.2 II.

\qed \par
 \enddemo

\proclaim {Lemma 4.2} \par
\it{Let $X:= \Bbb C^3/<G',T>$, and $f:\widetilde X \longrightarrow X$ the
crepant resolution as above. Then the Euler number of $\widetilde X$ is given
by
$$ \chi(\widetilde X) = \frac {|G'|-3(r-1)-1}{6}+3+2(r-1) \quad (|G'|=r^2)$$
\endproclaim

\demo{\bf{Proof}}

For an abelian group $G'$, we have a toric resolution
$$ \pi : \widetilde Y \longrightarrow Y = \Bbb C^3/G', $$
and $\chi(\widetilde Y) = \vert G'\vert$ ([8],[11]). \par
By the action of $\frak S_3$, the number of fixed points in the exceptional
divisor by $\sigma$ is equal to $r$, hence
$$ \chi(\widetilde Y/\frak S_3)=\frac {\vert G' \vert -3(r-1)-1}{6}+r. $$
By the resolution of the fixed loci, Eular characteristic of the exceptional
locus is 2 or 3. (Claim IV of \S 2 and Claim III of \S 4) \par
Therefore,
$$ \chi(\widetilde X) = \frac {\vert G'\vert -3(r-1)-1}{6} + 2(r-1)+3.  \qed$$
\enddemo

\proclaim{Theorem 4.3} \par
\it{$$\chi (\widetilde X) = \sharp \{ \text{conjugacy class of $G$} \}. $$
}\endproclaim
\demo{\bf{Proof}}

For 3(r-1) elements of type $\frac 1{r}(k_i,h_i,h_i)=:a_i,$ there are three
conjugate elements $a_i$, $Ta_iT^{-1}$ and $T^{-1}a_iT$. For other nontrivial
elements, there are six conjugate elements $a_i,$  $Ta_iT^{-1}$, $T^{-1}a_iT$,
$Sa_iS$, $STa_iT^{-1}S$ and $ST^{-1}a_iTS.$ There are $r-1$ conjugacy classes
of type $[a_iS]$ . There are 3 conjugacy classes $[e]$, $[S]$, $[T]$. So
the number of the conjugacy classes are: \par
$$ (r-1)+\frac 1{6}\{r^2-3(r-1)-1\}+ (r-1) +3. $$
Then
$$
   \chi(\widetilde X) = \sharp \{ \text{ conjugacy class of } G \}. \qed$$
\enddemo
\subhead \nofrills {\bf \S 5. Proof of IV} \endsubhead \par
\proclaim{Proposition 5.1}\par
\it{Let $X = {\Bbb C} ^3/ G$, and $Y = {\Bbb C}^3 /G'$. Then there exists the
following diagram: }
$$ \CD
   @. @. {\widetilde X} \\
   @. @. @VV{\tau}V     \\
   @. {\widetilde Y} @>{\widetilde\mu}>> {\widetilde Y}/{\frak S_3} \\
   @. @VV{\pi}V @VV{\widetilde\pi}V \\
   {\Bbb C^3} @>>> {\Bbb C^3}/{G'}=Y @>{\mu}>> {\Bbb C^3}/G=X
 \endCD $$
\it{where $\pi$ is a resolution of the singularity of $Y$, and $\widetilde\pi$
is the induced morphism, $\tau$ is a resolution of  the singularity by $\frak
S_3$, and $\tau \circ \widetilde \pi$ is a crepant resolution of the
singularity of $X$.}\endproclaim

\proclaim{Claim I} \it{There exists a toric resolution of $Y$ where $\frak S_3$
acts symmetrically on the exceptional divisors.}
\endproclaim
\demo{\bf{Proof}}
We construct a toric resolution as in \S 4.
\qed\enddemo
In the center of the triangle whose vertices are (1,0,0),(0,1,1) and (0,0,1),
there exist one $\Bbb P^2$ as an exceptional component. Then it is sufficient
to show the case $G'= < \frac 1{3}(1,1,1) >$.
\proclaim{Claim II}
\it{Let $F$ be the fixed locus on $\widetilde Y$ under the action of $\frak
S_3$, and $E$ be exceptional divisors of $\widetilde Y \longrightarrow Y$. Then
$$ F_0 := F \cap E =\text{3 points}. $$
\endproclaim
\demo {\bf{Proof}}
Considering the dual graph of exceptional divisors by the toric resolution,
 there is a possibility of the central component;
  a divisor which is isomorphic to $\Bbb P^2$.
Then the number of the fixed points under the action of $\frak S_3$ is three as
in \S 4.
\qed \enddemo
Futhermore the $\frak S_3$-action in the neighbourhood of a fixed point is
analytically isomorphic to some linear action. \par

\proclaim {Lemma 5.2} \par
\it{Let $X:= \Bbb C^3/<G',T>$, and $f:\widetilde X \longrightarrow X$ the
crepant resolution as above. Then the Euler number of $\widetilde X$ is given
by
$$ \chi(\widetilde X) = 9 $$
\endproclaim

\demo{\bf{Proof}}
For an abelian group $G'$, we have a toric resolution
$$ \pi : \widetilde Y \longrightarrow Y = \Bbb C^3/G', $$
and $\chi(\widetilde Y) = \vert G'\vert =3$ ([8],[11]). \par
By the action of $\frak S_3$, the number of fixed points in the exceptional
divisor by $\sigma$ is three, hence
$$ \chi(\widetilde Y/\frak A_3)=\frac 1{6} (\vert G' \vert -3)+3=3. $$
By the resolution of the fixed loci, Eular characteristic of the each
exceptional locus is 3. (Claim III in \S 4) \par
Therefore,
$$ \chi(\widetilde X) = \frac 1{6}(\vert G'\vert -3) + 3\times3=9.  \qed$$
\enddemo

\proclaim{Theorem 5.3} \par
\it{$$\chi (\widetilde X) = \sharp \{ \text{conjugacy class of $G$} \}. $$
}\endproclaim
\demo{\bf{Proof}}
There are nine conjugacy classes in $G$:\par
 identity, $\frac 1{3}(1,1,1)$,$ \frac 1{3}(2,2,2),$ $ S$,$ \frac
1{3}(1,1,1)S$,$ \frac 1{3}(2,2,2)S, T$,$ \frac 1{3}(1,1,1) $and
$\frac 1{3}(2,2,2).$ \qed \enddemo

Before we think about general case, we see the case of V in Theorem 1.2. \par

\proclaim{Theorem 5.4} \par
\it{For $G=G_5$, the conjecture II holds.}\endproclaim
\demo{\bf{Proof}}
Similarly we construct a crepant resolution:\par

$$ \CD
   @. @. {\widetilde X} \\
   @. @. @VV{\tau}V     \\
   @. {\widetilde Y} @>{\widetilde\mu}>> {\widetilde Y}/{\Bbb Z_2} \\
   @. @VV{\pi}V @VV{\widetilde\pi}V \\
   {\Bbb C^3} @>>> {\Bbb C^3}/{<C>}=Y @>{\mu}>> {\Bbb C^3}/G_5=X
 \endCD $$
There is an exceptional divisor which is isomorphic to $\Bbb P^2$ in
$\widetilde Y$. And there exist three singularities in $\widetilde Y/\Bbb Z_2$,
they become three $\Bbb P^1$-bundles in $\widetilde X$.
Then the Eular number of $\widetilde X$ is 6. \par
On the other hands, the number of the conjugacy classes in $G$ is 6: id, $C$,
$C^2$, $S$, $CS$ and  $C^2S$. \qed \enddemo

In general, we get the result:
\proclaim{Theorem 5.5}\par
\it{ As the normal subgroup $G'$ of $G$, we take an abelian group generated by
all of the diagonal matrices in $G$. Then
$$ \chi (\widetilde X) = \frac 1{2}(r^2-3r+2)+6r+3 ,$$
and it equals to the number of the conjugacy classes in $G$.
}\endproclaim
\demo{\bf{Proof}}
$$ \chi(\widetilde Y)= 3r^2. $$
There are 3+9(r-1) fixed points on the exceptional divisors in $\widetilde Y$,
and  they turn to the 3+3(r-1) singularities in $\widetilde Y/\frak S_3$.
$$ \chi(\widetilde Y/\frak S_3)=\frac 1{6}\{3r^2-9(r-1)-3\}+3(r-1)+3 $$
Then
$$ \chi(\widetilde X)=1{6}\{3r^2-9(r-1)-3\}+2\times3(r-1)+3\times3.$$
And this number coincide with the number of conjugacy classes, because
$G_4=G_3\amalg G_3C\amalg G_3C^2$. \qed
\enddemo

Therefore, Main theorem (Theorem 1.2) is proved!

\enddocument